\documentclass[12pt]{article}
\usepackage{fullpage,amssymb}
\usepackage{graphicx,psfrag}
\newcommand{\CC}{\mathbb{C}}
\newcommand{\RR}{\mathbb{R}}
\newcommand{\ZZ}{\mathbb{Z}}
\newcommand{\calm}{\mathcal{M}}

\begin{document}

\title{Goldman flows on the Jacobian}
\author{Lisa C. Jeffrey and David B. Klein \\
Mathematics Department, University of Toronto \\
Toronto, Ontario, Canada M5S 2E4 \\
jeffrey@math.toronto.edu, dklein@math.toronto.edu}
\date{}
\maketitle

\begin{abstract}
We show that the Goldman flows preserve the holomorphic structure
on the moduli space of homomorphisms of the fundamental group of
a Riemann surface into $U(1)$, in other words the Jacobian.
\end{abstract}




\section{Introduction}

This note concerns  the moduli space $\calm(G)$ of
conjugacy classes of homomorphisms of the fundamental group
of a compact orientable 2-manifold $\Sigma$ into a Lie group $G$.

This object has recently attracted a great deal of interest in symplectic
and algebraic geometry and mathematical physics.  In mathematical
physics it appears as the space of gauge equivalence classes of flat
connections on a 2-manifold.  In algebraic geometry it appears as the
moduli space of holomorphic bundles on a Riemann surface.

The smooth locus of the space $\calm(G)$ has a symplectic structure;
see, for instance, \cite{Goldman84} or \cite{Karshon}.
If the 2-manifold is equipped with a complex structure,
the space $\calm(G)$ inherits a  complex structure
compatible with the symplectic structure.
When $G=U(1)$, the moduli space of a Riemann surface coincides
with its Jacobian, which is a complex torus whose complex dimension
is the genus of the surface.

In \cite{Goldman86}, W. Goldman studied the Hamiltonian flows
of certain natural functions on the moduli space.
These functions are constructed from functions that send a flat
connection to its holonomy  around a specific simple closed
curve $C$ in the 2-manifold.

In this paper we show that when the gauge group is $U(1)$
the Goldman flows preserve the complex structure on $\calm(U(1))$.
When the gauge group is $SU(2)$, the Goldman flows
are ill defined on a set of real codimension 3,
(see \cite{JeffreyWeitsman}),
which is inconsistent with these flows
preserving the complex structure.

{\begin{figure}[t]
\psfrag{homology}{
 $\lambda_1,\dots,\lambda_{2g}
 \in H_1(\Sigma,\mathbb{Z})$}
\psfrag{fundamental}{
 $\widetilde{\lambda}_1,\dots,\widetilde{\lambda}_{2g}
 \in\pi_1(\Sigma)$}
\psfrag{L1}{${\scriptstyle \lambda_1}$}
\psfrag{L2}{${\scriptstyle \lambda_2}$}
\psfrag{CL2}{${\scriptstyle C=\lambda_2}$}
\psfrag{L2g1}{${\scriptstyle \lambda_{2g-1}}$}
\psfrag{L2g}{${\scriptstyle \lambda_{2g}}$}
\psfrag{0L1}{${\scriptstyle \widetilde{\lambda}_1}$}
\psfrag{0L2}{${\scriptstyle \widetilde{\lambda}_2}$}
\psfrag{C0L2}{${\scriptstyle C=\widetilde{\lambda}_2}$}
\psfrag{0L2g1}{${\scriptstyle \widetilde{\lambda}_{2g-1}}$}
\psfrag{0L2g}{${\scriptstyle \widetilde{\lambda}_{2g}}$}
\psfrag{C}{${\scriptstyle C}$}
\includegraphics{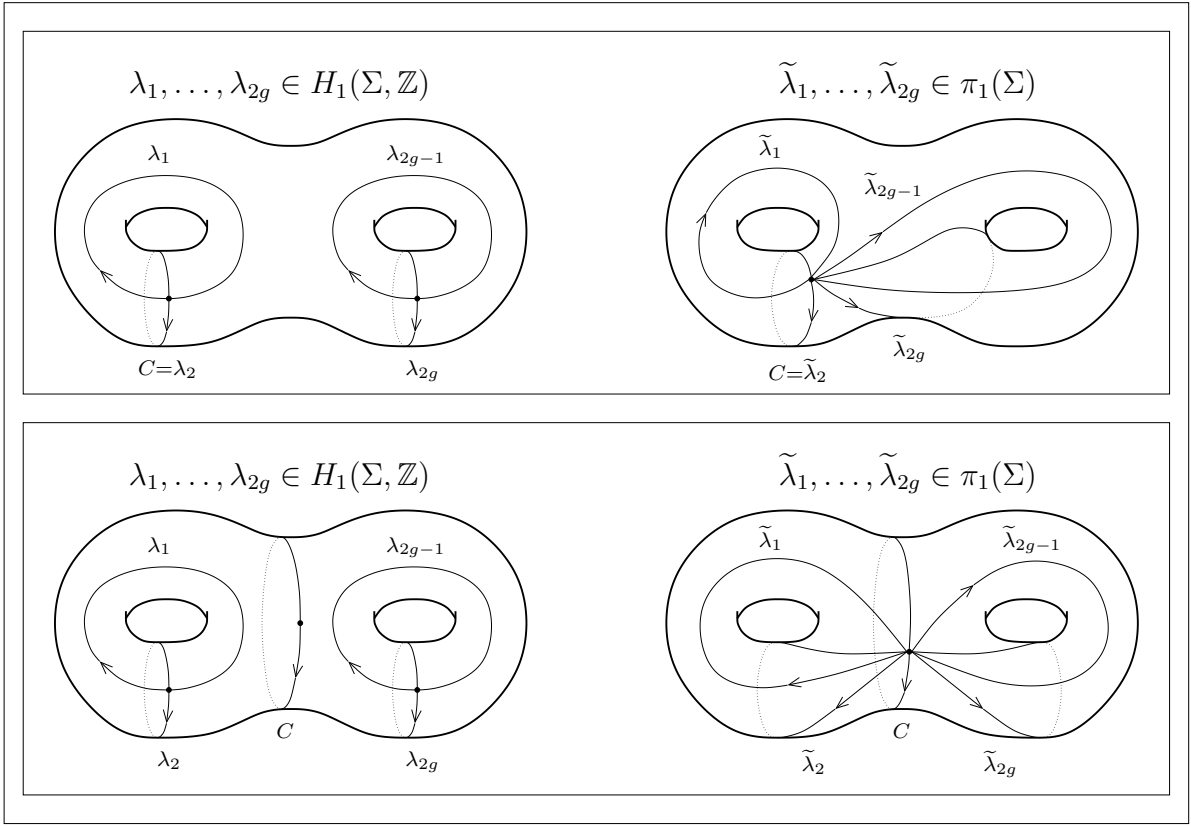}
\caption{\label{x_surfaces}
The curve $C$ is either nonseparating (top), or separating (bottom).}
\end{figure}}

\section{The $U(1)$ Goldman flow}

Let $G=U(1)$.  The Goldman flow on $\calm(G)$
is a periodic $\RR$-action $\{\Xi_s\}_{s\in\RR}$
associated to a simple closed curve $C$.  Since the Lie group
is abelian, the moduli space is ${\rm Hom}(\pi, G)$.

Choose a symplectic basis of $H_1(\Sigma,\ZZ)$, in other words
a collection of cycles $\{\lambda_1,\dots\lambda_{2g}\}$
in which all the intersections are empty except for
$\lambda_{2j-1}$ and $\lambda_{2j}$, which intersect
once transversely with positive intersection index.
If the curve $C$ is nonseparating then we let $\lambda_2=C$ in this
symplectic basis; if the curve $C$ is separating then we assume that
the cycles $\lambda_1,\dots,\lambda_{2g}$ do not intersect $C$.
Choose a basepoint on $C$ for the fundamental group of $\Sigma$,
and lift the cycles in the symplectic basis to loops
$\widetilde{\lambda}_1,\dots,\widetilde{\lambda}_{2g}\in\pi_1(\Sigma)$
that only intersect $C$ at their endpoints.
See Figure~\ref{x_surfaces}.

For $G=U(1)$, identify $\Phi\in{\rm Hom}(\pi_1(\Sigma),G)$ with
$(\Phi_1,\dots,\Phi_{2g})\in G^{2g}$ by letting
$\Phi_j=\Phi(\widetilde{\lambda}_j)$.
If the simple closed curve $C$ is nonseparating then
the associated Goldman flow on ${\rm Hom}(\pi, G)=G^{2g}$ is
$$ \Xi_s(\Phi)=(e^{2\pi is}\Phi_1,\Phi_2,\dots,\Phi_{2g}), $$
for $s\in\RR/\ZZ$.
If $C$ is separating then the Goldman flow is trivial,
$\Xi_s(\Phi)=\Phi$.

The Goldman flows can be described using gauge theory
as follows, (cf. \cite{Goldman04}).
If $A$ is a flat connection on the trivial $G$-bundle
over $\Sigma$, and $A_j$ is the holonomy of $A$ along
the loop $\widetilde{\lambda}_j$, then the map that sends $A$ to
$(A_1,\dots,A_{2g})\in G^{2g}={\rm Hom}(\pi_1(\Sigma),G)$
identifies the space of gauge equivalence classes of
flat connections with the moduli space; see \cite{Jeffrey}.
Let $\widehat{\Sigma}$ be the complement of $C$ in $\Sigma$,
and let $U_- \cup U_+$ be the intersection of $\widehat{\Sigma}$
with a tubular neighbourhood of $C$ in $\Sigma$.  The open sets
$U_-$ and $U_+$ can be thought of as neighbourhoods of the two
``boundary components" of $\widehat{\Sigma}$, as shown in
Figure~\ref{x_gauge}.  For a flat connection $A$ on $\Sigma$, define
$$ \Xi_s(A)=A^{g_s}, $$
where $g_s$ is a gauge transformation on $\widehat{\Sigma}$
with $g_s=1$ on $U_+$ and $g_s=e^{2\pi is}$ on $U_-$.
Here, since $g_s$ is constant on both $U_+$ and $U_-$
and since the gauge group is abelian,
the flat connection $A^{g_s}$ on $\widehat{\Sigma}$
extends (by $A$) to a flat connection $\Xi_s(A)$ on $\Sigma$.
If the curve $C$ separates $\Sigma$ into two components
then $\Xi_s(A)=A$ because $g_s$ may be chosen to be
a locally constant gauge transformation on $\widehat{\Sigma}$,
which acts trivially since $G$ is abelian.
If the curve $C$ is nonseparating, however, then the
gauge transformations $g_s$ on $\widehat{\Sigma}$
act nontrivially and do not come from gauge transformations
on $\Sigma$ for $s\notin\ZZ$, so in general the $\Xi_s(A)$
are distinct elements of the moduli space $\calm(G)$,
(although they are the same when viewed as
elements of the moduli space of $\widehat{\Sigma}$).

{\begin{figure}[t]
\psfrag{connected}{$\widehat{\Sigma}$ connected}
\psfrag{disconnected}{$\widehat{\Sigma}$ disconnected}
\psfrag{U-}{${\scriptstyle U_-}$}
\psfrag{U+}{${\scriptstyle U_+}$}
\psfrag{g-}{${\scriptstyle g_s=e^{2\pi is}}$}
\psfrag{g+}{${\scriptstyle g_s=1}$}
\includegraphics{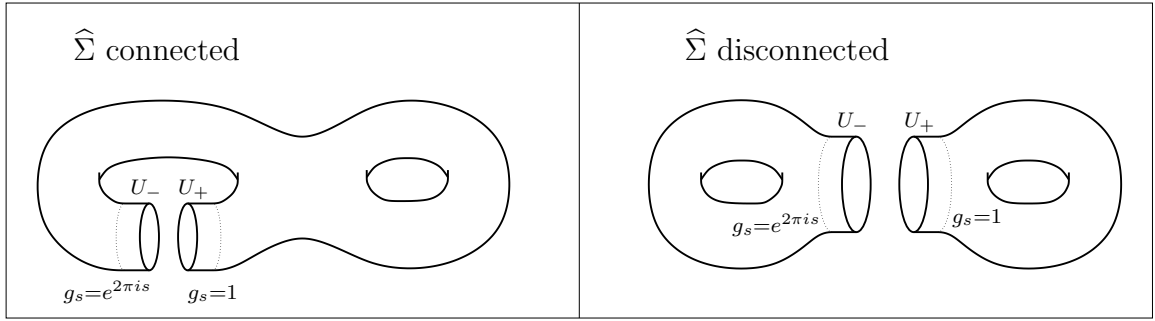}
\caption{\label{x_gauge}
 The gauge transformation $g_s$.}
\end{figure}}

\section{The Jacobian}

When $G=U(1)$, the moduli space $\calm(U(1))$ is
the Jacobian ${\rm Jac}(\Sigma)\cong U(1)^{2g}$.
The Jacobian inherits a complex structure from
the Riemann surface $\Sigma$, and identifies with
$\CC^g/\Lambda$ as a complex manifold
for a lattice $\Lambda$ described below.
See, for example, \cite{GriffithsHarris}.

The Jacobian is defined as
$$ {\rm Jac}(\Sigma)=\frac{H^0(\Sigma, K)^*}{H_1(\Sigma,\ZZ)}, $$
where $H_1(\Sigma,\ZZ)$ maps to $H^0(\Sigma, K)^*$ by integration:
a class $\lambda\in H_1(\Sigma,\ZZ)$
sends $\omega\in H^0(\Sigma,K)$ to $\int_\lambda\omega\in\CC$.
Explicitly, choose a basis $\{\omega_1,\dots,\omega_g\}$
of $H^0(\Sigma,K)$, and use the dual basis to identify
$H^0(\Sigma,K)^*$ with $\CC^g$.
Let $F$ be the resulting map from $H_1(\Sigma,\ZZ)$ to $\CC^g$,
$$ F(\lambda)=\left(\begin{array}{c} \int_\lambda\omega_1 \\
\vdots \\ \int_\lambda\omega_g \end{array}\right), $$
and equate $H_1(\Sigma,\ZZ)\subset H^0(\Sigma,K)^*$ with the lattice
$$ \Lambda=\{F(\lambda) \ : \
\lambda\in H_1(\Sigma,\ZZ)\}\subset\CC^g. $$

A choice of basis $\{\lambda_1,\dots,\lambda_{2g}\}$ of
$H_1(\Sigma,\ZZ)$ identifies $\CC^g/\Lambda$ with
$U(1)^{2g}$ as follows.
Viewed as a \emph{real} vector space, $\CC^g$ is
spanned by $\{F(\lambda_1),\dots,F(\lambda_{2g})\}$,
$$ \CC^g=\left\{\sum_{j=1}^{2g}v_j F(\lambda_j) \ : \
v_j\in\RR\right\}. $$
Identify $\CC^g/\Lambda$ with $U(1)^{2g}$
by the group isomorphism that maps
$$ \left[\sum_{j=1}^{2g}v_j F(\lambda_j)\right]\in\CC^g/\Lambda $$
to
$$ (\exp 2\pi iv_1,\dots,\exp 2\pi iv_{2g}). $$
We define $z_j=\exp 2\pi iv_j$.


\section{Goldman flows on the Jacobian}

The Goldman flow on the Jacobian is defined as follows.
If $C$ is a nonseparating simple closed curve then
choose $C$ as the generator $\lambda_2$ in a
symplectic basis $\{\lambda_1,\dots,\lambda_{2g}\}$ of
$H_1(\Sigma,\ZZ)$.  So the Goldman flow associated to $C$ is
$$ (z_1,z_2,\dots,z_{2g})\mapsto(e^{2\pi is}z_1,z_2,\dots,z_{2g}) $$
for $s\in\RR/\ZZ$.
This corresponds to translation by $sF(\lambda_1)$ in $\CC^g$,
which is clearly a holomorphic self map of $\CC^g/\Lambda$.

If $C$ is a separating simple closed curve then the Goldman flow
on the Jacobian is trivial (because the gauge group is abelian).


\end{document}